\documentclass[11pt,twoside]{article}
\usepackage[left=3cm,right=3cm,top=3cm,bottom=3cm]{geometry}
\usepackage{amssymb}
\usepackage{graphicx}
\usepackage{graphicx,type1cm,eso-pic,color}
\usepackage{amssymb}
\usepackage{amssymb,amsmath}
\usepackage{graphicx,epsfig}
\usepackage{enumerate}
\usepackage{color}
\usepackage {multicol}
\usepackage{lineno}
\numberwithin{equation}{section}
\newtheorem{theorem}{Theorem}

\newtheorem{example}{Example}

\makeatletter
\AddToShipoutPicture{
	\setlength{\@tempdimb}{.5\paperwidth}
	\setlength{\@tempdimc}{.5\paperheight}
	\setlength{\unitlength}{1pt}
	\put(\strip@pt\@tempdimb,\strip@pt\@tempdimc)
}
\makeatother

\usepackage{scalerel,stackengine}
\stackMath
\newcommand\reallywidehat[1]{%
	\savestack{\tmpbox}{\stretchto{%
			\scaleto{%
				\scalerel [\widthof{\ensuremath{#1}}]{\kern-.6pt\bigwedge\kern-.6pt}%
				{\rule[-\textheight/2]{1ex}{\textheight}}%WIDTH-LIMITED BIG WEDGE
			}{\textheight}% 
		}{0.5ex}}%
	\stackon[1pt]{#1}{\tmpbox}%
}
\begin{document}
	\setcounter{page}{1}
	%\newcounter{equation}[section]
	\thispagestyle{empty}
	\markboth{}{}

	\pagestyle{myheadings}
	\markboth{}{ }
	
	\date{}
	
	%\maketitle
	
	\noindent  
	
	\vspace{.1in}
	
	{\baselineskip 20truept
		
		\begin{center}
			{\Large {\bf Extension of Yager's negation of a probability distribution based on uncertainty measures }} \footnote{\noindent
                 {\bf $^{1}$} E-mail: skchaudhary1994@gmail.com
				{\bf $^{2}$} E-mail: nitin.gupta@maths.iitkgp.ac.in\\
				{\bf $^{2^{*}}$} corresponding author E-mail: pradeep.maths@kgpian.iitkgp.ac.in}\\
	\end{center}}
	\vspace{.1in}
	
	\begin{center}
		{\large {\bf Santosh Kumar Chaudhary$^1$}}
        {\large {\bf Pradeep Kumar Sahu$^{2^{*}}$ and Nitin Gupta$^{2}$}}\\
		\vspace{0.1cm}
		{\large {\it $^{1}$ Department of Statistics, Central University of Jharkhand, Cheri-Manatu, Ranchi, Jharkhand, 835222, India. }}\\
        {\large {\it $^{2,2^{*}}$ Department of Mathematics, Indian Institute of Technology Kharagpur, West Bengal 721302, India. }}\\
		\end{center}
	
	\vspace{.1in}
	\baselineskip 12truept
	%\begin{abstract}

	%\end{abstract}
	%\vspace{.1in}
	%\noindent  {\bf Key Words}: {\it residual entropy, past entropy, Varma entropy, Tsallis entropy.}
	\begin{center}
		{\bf \large Abstract}\\
	\end{center}
	 Existing research on negations primarily focuses on entropy and extropy. Recently, new functions such as varentropy and varextropy have been developed, which can be considered as extensions of entropy and extropy. However, the impact of negation on these extended measures, particularly varentropy and varextropy, has not been extensively explored. To address this gap, this paper investigates the effect of negation on Shannon entropy, varentropy, and varextropy. We explore how the negation of a probability distribution influences these measures, showing that the negated distribution consistently leads to higher values of Shannon entropy, varentropy, and varextropy compared to the original distribution. Additionally, we prove that the negation of a probability distribution maximizes these measures during the process. The paper provides theoretical proofs and a detailed analysis of the behaviour of these measures, contributing to a better understanding of the interplay between probability distributions, negation, and information-theoretic quantities.
	\section{Introduction}
   Negation is a pivotal method in information processing, especially when handling uncertain or incomplete data. Typically, when faced with such information, individuals tend to focus more on the positive aspects, often overlooking the significance of the negative components. However, in many cases, analyzing a problem from a negative perspective can yield valuable insights that might be more challenging to uncover by examining only the positive side. In fact, there are scenarios where information that is difficult to interpret from a positive viewpoint can be more easily extracted from a negative one. Furthermore, a more thorough and accurate understanding of a problem can be achieved by considering both the positive and negative aspects together, as this approach allows for a more comprehensive and well-rounded view of the situation.

In this context, we wish to address the challenge of determining the negation of a probability distribution. Specifically, we will explore how negation affects the key measures associated with probability distributions. In particular, our focus will be on varextropy and varentropy, which are important metrics for understanding the changes in the distribution's characteristics when it is negated. These measures will provide insight into how negation can influence the uncertainty represented by the distribution and how it can be used to refine the analysis of uncertain data.

Negation methods are extensively utilized in evidence theory, a framework designed to handle reasoning and decision-making under conditions of uncertainty. In evidence theory, negation techniques are employed to adjust models that are based on uncertain or incomplete information, helping to better capture and understand the underlying uncertainty. Among the many methods available for negation, entropy plays a particularly crucial role. For instance, Yager (2015) introduced a technique for negating probability distributions based on  GINI entropy and explored Tsallis entropy in this context. These entropy-based methods provide a systematic way to assess and process uncertainty, offering valuable insights into how negating a probability distribution alters the level of uncertainty represented by the data.

In addition to entropy, a concept that has recently gained increasing attention is extropy. Extropy is often viewed as the dual counterpart to entropy and serves a similar function in measuring uncertainty. According to the interpretation provided by Lad, Sanfilippo, and Agro (2015), the probability distribution that maximizes extropy is a uniform distribution. While extropy is mathematically related to entropy, it is expressed differently, offering an alternative approach to quantifying uncertainty. Given the numerous advantages of extropy, it has become an important area of exploration, particularly in the context of negating information. The ability to use extropy to analyze the effects of negation on probability distributions opens up new avenues for handling uncertainty in a more nuanced way.
Varentropy is a metric used to quantify the variability in the information content of a random vector and is unaffected by affine transformations.

\subsection{Shannon Entropy}
The entropy of a discrete probability distribution $ P = \{ p_1, \dots, p_n \} $, is defined as (Shannon 1948),

\[
H(P) = - \sum_{i=1}^{n} p_i \ln p_i
\]

Let $ X $ be an absolutely continuous random variable with probability density function (pdf) $ f(x) $. Define $ l_X = \inf \{ x \in \mathbb{R} : F(x) > 0 \} $, $ u_X = \sup \{ x \in \mathbb{R} : F(x) < 1 \} $, and $ S_X = (l_X, u_X) $. 
\subsection{Varentropy}
The varentropy of a discrete probability distribution $ P = \{ p_1, \dots, p_n \} $, is defined as (see, Bobkov
and Madiman (2011), Kontoyiannis and Verdu (2014), Arikan (2016), Di Crescenzo and Paolillo (2021), and Maadani et al. (2021))
      \[VH(P)= \sum_{i=1}^{n} p_i \left( \ln(p_i) \right)^2 -\left(  \sum_{i=1}^{n} p_i \ln(p_i) \right)^2 \] 
 
Varentropy serves as a measure of the variability in the information content.

\subsection{Extropy}

Lad et al. (2015) introduced a concept called extropy, which is the complement of Shannon entropy. The extropy of a discrete probability distribution $ P = \{ p_1, \dots, p_n \} $, is defined as (Lad, Sanfilippo, and Agro, 2015),

\[
J(P) = -\sum_{i=1}^{n} (1 - p_i) \ln(1 - p_i)
\]

 \subsection{Varextropy}
The varextropy of a discrete probability distribution $ P = \{ p_1, \dots, p_n \} $, is defined as (see Vaselabadi et al. (2021), Goodarzi (2022)
and Zaid et al. (2022))

      \[VJ(P)= \sum_{i=1}^{n} (1-p_i) \left( \ln((1-p_i)) \right)^2 -\left(  \sum_{i=1}^{n} (1-p_i) \ln((1-p_i)) \right)^2 \]

Varextropy also serves as a measure of the variability in the information content.

In this paper, we study negating a probability distribution that is based on  varextropy and varentropy . This approach provides a fresh perspective on how negation can be applied to probability distributions, offering new tools for analyzing and understanding uncertainty. By leveraging these advanced concepts, this method allows for a deeper exploration of how information can be negated and how these changes in negation affect the overall distribution, ultimately enhancing our ability to process uncertain data more effectively.

\section{Negation of a discrete probability distribution}
Yager (2015) investigated the negation of a given probability distribution using the entropy function. Yager (2015) explored how the knowledge within the negation of a probability distribution can be represented. He proposed a transformation method to derive the negation of a probability distribution and examined its properties. By applying the Dempster–Shafer theory of evidence, Yager (2015) demonstrated that among all possible negations, the one suggested in his work exhibits the maximal type of entropy.
While there are many different measures of entropy, Yager (2015) used the following to measure the entropy of a probability distribution,

\[
H_1(P) =  \sum_{i=1}^{n} (1 - p_i) p_i = 1 - \sum_{i=1}^{n} p_i^2.
\]
Yager (2015) selected this form of entropy measure instead of the classic Shannon entropy due to the simplicity of calculation it offers, as it does not involve logarithms.

Liu and Xiao (2024) investigated the negation of a given probability distribution using the extropy function. They demonstrated that through repeated negation of the probability, both the probability distribution and its associated extropy converge to stable values.

\subsection{Negation method}
Yager (2015) considered the problem of finding the
negation of a probability distribution and suggested following a negation of a discrete probability distribution. 
Assume $ X = \{ x_1, \dots, x_n \} $ and $ P = \{ p_1, \dots, p_n \} $ is the probability  distribution of $ X $ such that $ \sum_{i=1}^{n} p_i = 1 $, and $ p_i \in [0, 1] $. Let $ \bar{P} = \{ \bar{p}_1, \bar{p}_2, \dots, \bar{p}_n \} $ represent the negation of the probability distribution $ P $. The negation of the probability is given by (Yager, 2015) $
\bar{p}_i = \frac{1 - p_i}{n - 1}.$
$\bar{P}$ is also a probability distribution (Yager, 2015) since $
\sum_{i=1}^{n} \bar{p}_i =  1 \text{ and } \quad \bar{p}_i \in [0, 1].$
The inverse is irreversible in general, which means $ p_i \neq \bar{\bar{p}}_i $ but $ p_i = \bar{\bar{p}}_i $ when $n=2$ (Yager 2015). Moreover, $0 \leq \bar{p}_i \leq \frac{1}{n-1}, \quad i = 1, 2, \dots, n$  (Liu and Xiao (2024)). The fact that the operation is irreversible in general adds an interesting aspect to the analysis of uncertainty. The negation of probability distributions may have applications in Risk Analysis, Machine Learning, Data Science, Information Theory and Decision Theory.

\subsection{Negation with Shannon Entropy}
  The negation of entropy 
\[
H(P) = -\sum_{i=1}^{n} p_i \ln p_i
\]
is given as 

\[
H(\bar{P}) =- \sum_{i=1}^{n} \bar{p}_i \ln \bar{p}_i
\]

\subsection{Negation with varentropy}
The negation of varentropy
\[VH(P)= \sum_{i=1}^{n} p_i \left( \ln(p_i) \right)^2 -\left(  \sum_{i=1}^{n} p_i \ln(p_i) \right)^2 \] 
is given as 
  \[VH(\bar{P})= \sum_{i=1}^{n} \bar{p}_i \left( \ln(\bar{p}_i) \right)^2 -\left(  \sum_{i=1}^{n} \bar{p}_i \ln(\bar{p}_i) \right)^2 \]

\subsection{Negation with Varextropy}
The negation of varextropy
\[VJ(P)= \sum_{i=1}^{n} (1-p_i) \left( \ln((1-p_i)) \right)^2 -\left(  \sum_{i=1}^{n} (1-p_i) \ln((1-p_i)) \right)^2 \]
is given as 
      \[VJ(\bar{P})= \sum_{i=1}^{n} (1-\bar{p}_i) \left( \ln((1-\bar{p}_i)) \right)^2 -\left(  \sum_{i=1}^{n} (1-\bar{p}_i) \ln((1-\bar{p}_i)) \right)^2 \] 
\section{Examples}\label{sec_example}

\begin{example}
Let $ n = 2 $, $ X = \{x_1, x_2\} $, and $ P = \{p_1, p_2\} $. Let $ p_1 = p_2 = \frac{1}{2} $, then $
    H(P) =H(\bar{P})=  \ln 2, VH(P)=VH(\bar{P})  = 0$ and $VJ(P) =VJ(\bar{P})= 0.$
\end{example}

\begin{example}
Let $ n = 2 $, $ X = \{x_1, x_2\} $, and $ P = \{p_1, p_2\} $. Let $ p_1 = \frac{2}{5}, p_2 = \frac{3}{5} $, then
$H(P) = 0.470, \  H(\bar{P}) = 0.470, \ VH(P) = 0.0857, \ VH(\bar{P}) = 0.0857, \   VJ(P) = 0.0857$ and $  VJ(\bar{P}) = 0.0857.$
  
\end{example}

\begin{example}
Let $ n = 2 $, $ X = \{x_1, x_2\} $, and $ P = \{p_1, p_2\} $. Let $ p_1 = \frac{1}{10}, p_2 = \frac{9}{10} $, then  $ H(P) =  0.325, \  H(\bar{P}) = 0.325, \    VH(P) =  0.0557, \    VH(\bar{P}) = 0.0557, \  VJ(P) =  0.0557 $   and $ VJ(\bar{P}) = 0.0557$
\end{example}
Since there are only two basic events, entropy, varentropy, and varextropy, they will not change after probability distribution negation iteration. As the number of iterations increases, the entropy, varentropy and varextropy remain unchanged for $n=2.$ That is, $H(P) =H(\bar{P})= H(\bar{\bar{P}})=\dots, \  VH(P) =VH(\bar{P})= VH(\bar{\bar{P}})=\dots$ and $VJ(P) =VJ(\bar{P})= VJ(\bar{\bar{P}})=\dots.$

\begin{figure}[htp]
    \centering
    \includegraphics[width=10cm]{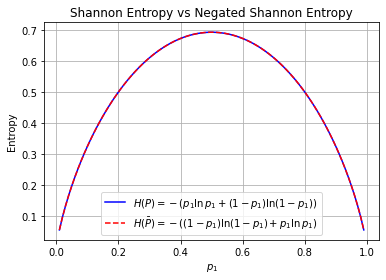}
    \caption{Change of $H(P)$ and $H(\bar{P})$}
    \label{fig:1a}
\end{figure}

\begin{figure}[htp]
    \centering
    \includegraphics[width=10cm]{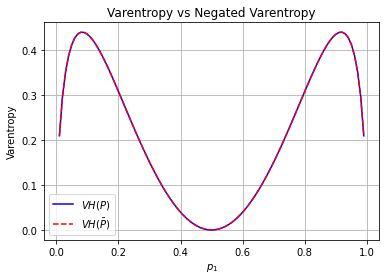}
    \caption{Change of $VH(P)$ and $VH(\bar{P})$}
    \label{fig:1b}
\end{figure}

\begin{figure}[htp]
    \centering
    \includegraphics[width=10cm]{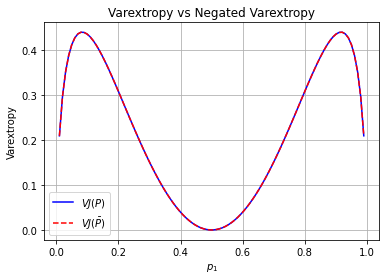}
    \caption{Change of $VJ(P)$ and $VJ(\bar{P})$}
    \label{fig:1c}
\end{figure}

Figure 1, Figure 2 and Figure 3 shows the change of Shanon entropy, varentropy and varextropy, respectively, with the change of $p_1$ from 0 to 1 for n=2. Since $H(\bar{P})=$H(P), therefore gragh coincides. The entropy will be highest when $p_1=p_2=0.5$ because the distribution is most uncertain. As $p_1$ approaches 0 or 1, the entropy decreases, reaching 0 when one of the probabilities is 0 (i.e., the distribution becomes deterministic).

\begin{example}
    For \( n \) equally likely outcomes (uniform distribution), each \( p_i \) is given by,

\[
p_i = \frac{1}{n} \quad \text{for all} \quad i = 1, 2, \dots, n; \  n>1.
\]

For the negated distribution \( \bar{P} \), the probabilities are,

\[
\bar{p}_i = \frac{1 - p_i}{n - 1} = \frac{1 - \frac{1}{n}}{n - 1} = \frac{n - 1}{n(n - 1)} = \frac{1}{n}
\]

Thus, the entropy for \( P \), which is a uniform distribution, is,

\[
H(P) = - \sum_{i=1}^{n} p_i \ln(p_i) = - n \cdot \frac{1}{n} \ln \left( \frac{1}{n} \right) = \ln(n)
\]

Since \( \bar{p}_i \) is also uniform, the entropy of \( \bar{P} \) is,

\[
H(\bar{P}) = \ln(n)
\]
Similarly, 
\[
VH(P) = \sum_{i=1}^{n} p_i \left( \ln(p_i) \right)^2 - \left( \sum_{i=1}^{n} p_i \ln(p_i) \right)^2= (\ln n)^2 - (- \ln n)^2 = (\ln n)^2 - (\ln n)^2 = 0
\]
and 
\[
VH(\bar{P}) = 0
\]
\end{example}
Both Shannon entropy and negated Shannon entropy curves will overlap because the entropies are identical for a uniform distribution. We observe thatas \( n \) increases, both \( H(P) \) and \( H(\bar{P}) \) increase. Figure 4 shows changes in \( H(P) \) and \( H(\bar{P}) \)  for different value of $n$. Figure 5 shows that \( H(P) \) = \( H(\bar{P}) \)=0 for all $n$.

\begin{example}
    The Varextropy of a discrete probability distribution \( P = \{ p_1, p_2, \dots, p_n \} \) is defined as,
\[
VJ(P) = \sum_{i=1}^{n} (1 - p_i) \left( \ln(1 - p_i) \right)^2 - \left( \sum_{i=1}^{n} (1 - p_i) \ln(1 - p_i) \right)^2
\]
For a uniform distribution, where \( p_i = \frac{1}{n} \) for all \( i = 1, 2, \dots, n \).
The varentropy for \( P \) is,
\[
VJ(P)= VJ(\bar{P})  = n \left( 1 - \frac{1}{n} \right) \left( \ln\left(1 - \frac{1}{n}\right) \right)^2 - \left( n \left( 1 - \frac{1}{n} \right) \ln\left(1 - \frac{1}{n}\right) \right)^2
\]
\end{example}
Figure 6 shows changes in \( VJ(P) \) and \( VJ(\bar{P}) \)  for different value of $n$.  \( VJ(P) \) is decreasing when n increases.

\begin{figure}[htp]
    \centering
    \includegraphics[width=10cm]{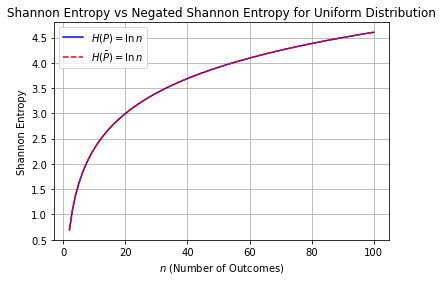}
    \caption{Change of $H(P)$ and $H(\bar{P})$ for different $n$ }
    \label{fig:2a}
\end{figure}

\begin{figure}[htp]
    \centering
    \includegraphics[width=10cm]{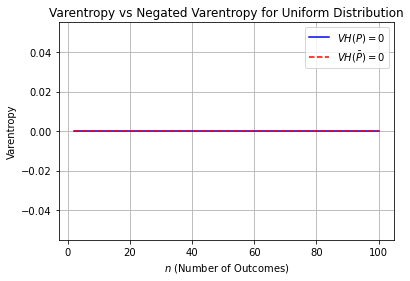}
    \caption{Change of $VH(P)$ and $VH(\bar{P})$ for different $n$ }
    \label{fig:2b}
\end{figure}

\begin{figure}[htp]
    \centering
    \includegraphics[width=10cm]{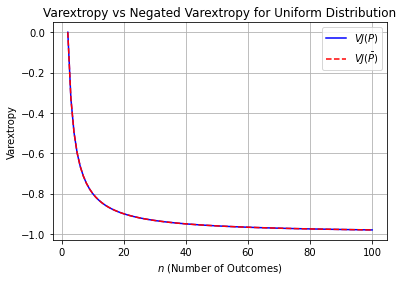}
    \caption{Change of $VJ(P)$ and $VJ(\bar{P})$ for different $n$ }
    \label{fig:2c}
\end{figure}

\begin{example}
Let $ n = 3 $, $ X = \{x_1, x_2, x_3\} $, and $ P = \{p_1, p_2, p_3\} $. Let $ p_1 = \frac{1}{3}, p_2 = \frac{1}{3}, \ p_3 = \frac{1}{3} $, then $    H(P) =H(\bar{P}) = \ln 3, \     VH(P) =   VH(\bar{P}) = 0, \  \text{and} \  VJ(P) =   VJ(\bar{P}) = 0.0837. $
\end{example}

\begin{example}
Let $ n = 3 $, $ X = \{x_1, x_2, x_3\} $, and $ P = \{p_1, p_2, p_3\} $. Let $ P = \left\{ \frac{1}{10}, \frac{3}{10}, \frac{6}{10} \right\} $, and $ \bar{P} = \left\{ \frac{9}{20}, \frac{7}{20}, \frac{1}{5} \right\} $. Then,
$H(P) = 0.898, \    H(\bar{P}) =  0.4236, \ VH(P) = 0.3146, \  VH(\bar{P}) = 0.3146, \ VJ(P) =  0.0837$ and $ VJ(\bar{P}) = 0.0837.$
\end{example}

\begin{example}\label{3a}
    Let the probability distribution be $ P = \{ p_1, p_2, p_3, p_4 \} = \{ 0.4, 0.3, 0.2, 0.1 \} $, and the negated distribution $ \bar{P} = \{ \bar{p}_1, \bar{p}_2, \bar{p}_3, \bar{p}_4 \} $, where
\[
\bar{p}_i = \frac{1 - p_i}{3} \quad \text{for} \quad i = 1, 2, 3, 4.
\]
Thus, the negated probabilities are,
\[
\bar{p}_1 = \frac{0.6}{3} = 0.2, \quad \bar{p}_2 = \frac{0.7}{3} = 0.2333, \quad \bar{p}_3 = \frac{0.8}{3} = 0.2667, \quad \bar{p}_4 = \frac{0.9}{3} = 0.3.
\]
Then,
 $H(P)=1.2799,\ J(X)=0.8295, VH(X)=0.1809, VJ(X)=-0.3926$

$ H(\bar{P}) = 1.3751, J(\bar{P})=0.8593, VH(\bar{P}) =0.0220, VJ(\bar{P})  =-0.4849 $

$ H(\bar{\bar{P}})= 1.3851,   J(\bar{\bar{P}})=0.8626, VH(\bar{\bar{P}})=0.0025,   VJ(\bar{\bar{P}})=-0.4953 $

 $H(\bar{\bar{\bar{P}}})=1.3862,  J(\bar{\bar{\bar{P}}})=0.8630,  VH(\bar{\bar{\bar{P}}})= 0.0003,   VJ(\bar{\bar{\bar{P}}})=-0.4964$

\begin{figure}[htp]
    \centering
    \includegraphics[width=12cm]{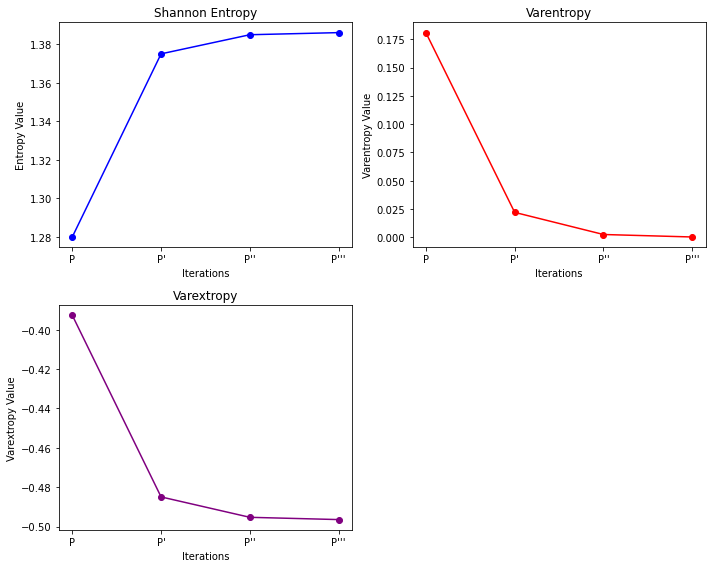}
    \caption{Plot of numerical value for Example \ref{3a} }
    \label{fig:3a}
\end{figure}

\end{example}

\begin{example}\label{3b}
Let $ X = \{ x_1, x_2, x_3 \} $ and the probability distribution $ P = \{ p_1, p_2, p_3 \} $, where $ p_1 = 0.6 $, $ p_2 = 0.3 $, and $ p_3 = 0.1 $. The negated distribution $ \bar{P} = \{ 0.2, 0.35, 0.45 \} $. Then,
$ H(P) = 1.1023,  H(\bar{P}) = 1.0486, \   VH(P) = -0.3937,    VH(\bar{P}) = 0.0917, \ VJ(P) = 0.118 $   and $VJ(\bar{P}) = 0.118.$

Then
 $H(P)=0.8979,\ VH(X)=0.3153 , \ VJ(X)= -0.0707$

$ H(\bar{P}) =1.0487 , \  VH(\bar{P}) = 0.0911, \ VJ(\bar{P}) =  -0.2629 $

$ H(\bar{\bar{P}})= 1.0868,   VH(\bar{\bar{P}})=0.0235,   VJ(\bar{\bar{P}})=-0.3121 $

 $H(\bar{\bar{\bar{P}}})=1.0956 ,   VH(\bar{\bar{\bar{P}}})= 0.0059,   VJ(\bar{\bar{\bar{P}}})= -0.3246 $
 
 $H(\bar{\bar{\bar{\bar{P}}}})= 1.0979,   VH(\bar{\bar{\bar{\bar{P}}}})= 0.0015,   VJ(\bar{\bar{\bar{\bar{P}}}})= -0.3278.$

\begin{figure}[htp]
    \centering
    \includegraphics[width=12cm]{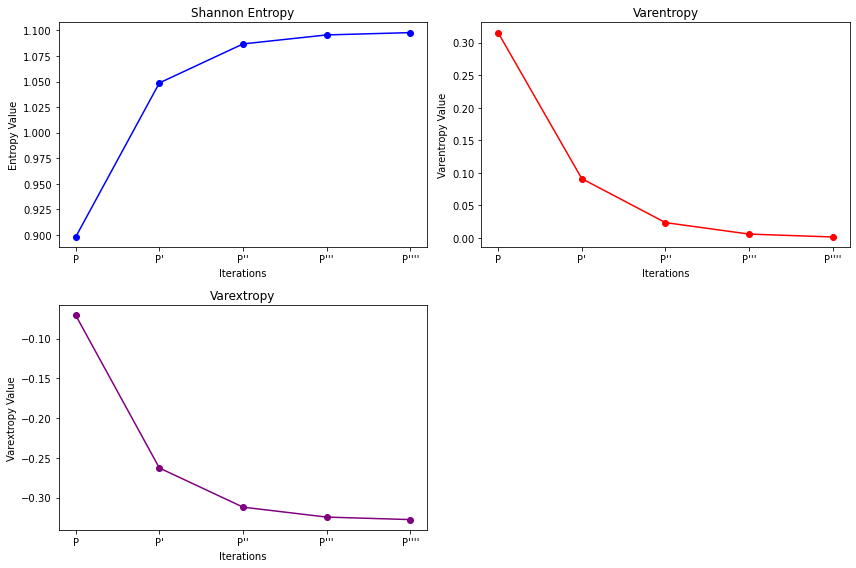}
    \caption{Plot of numerical value for Example \ref{3b} }
    \label{fig:3b}
\end{figure}

\end{example}
\begin{theorem}
 Assume the event space  $ X = \{ x_1, x_2, \dots x_n \} $ and the probability distribution $ P = \{ p_1, p_2, \dots, p_n \} $, $\bar{P}$ represents the inverse of P, then $H( \bar{P})\geq H(P).$
\end{theorem}
\noindent \textbf{Proof} We have seen in Section \ref{sec_example} that $H(P)=H(\bar{P})$  for $n=2.$ For $n\geq 3,$ the difference 
\begin{align}\label{lagrange1}
  Y_1=H(P)-H(\bar{P})  =& -\sum_{i=1}^{n} p_i \ln p_i+ \sum_{i=1}^{n} \bar{p}_i \ln \bar{p}_i \nonumber\\
 =&  -\sum_{i=1}^{n} p_i \ln p_i+ \sum_{i=1}^{n} \left(\frac{1-p_i}{n-1}\right) \ln \frac{1-p_i}{n-1} 
\end{align}
The Lagrange function under equation \ref{lagrange1} is
\begin{align}
T_1=  -\sum_{i=1}^{n} p_i \ln p_i+ \sum_{i=1}^{n} \left(\frac{1-p_i}{n-1}\right) \ln \frac{1-p_i}{n-1} + \lambda\left( \sum_{i=1}^{n} p_i -1 \right) . 
\end{align}
The partial derivative of $T_1$ with respect to $p_i$ is
\begin{align}
\frac{\partial T_1}{\partial p_i} = -\ln(p_i)-\frac{1}{n-1} \ln\left(\frac{1-p_i}{n-1}\right) - \frac{n}{n-1} +\lambda, \ \forall \  i=1,2, \dots, n. \label{3.3eqn}
\end{align}
and 
\begin{align}
\frac{\partial T_1}{\partial \lambda}= \sum_{i=1}^{n} p_i -1 
\end{align}
Lets solve  $\frac{\partial T_1}{\partial p_i}=0$ and $\frac{\partial T_1}{\partial \lambda}=0.$ to get stationary points. Equation \ref{3.3eqn} can be written as 
\begin{align}
    -\ln(p_i)-k_1 \ln\left(1-p_i\right) = k_2.
\end{align}
Let \begin{align}
    D_1=-\ln(p_i)-k_1 \ln\left(1-p_i\right),\\
    \frac{\partial D}{\partial p_i}=- \frac{1}{p_i}+k_1\frac{1}{1-p_i} <0.
\end{align}
Since $D_1$ decreases in a monotonic manner, the solution to $D = k_2$ is unique, implies that $p_i = 1/n$  for $i = 1, 2, ..., n.$

Hence, $T_1$ has maximum value 0 when $p_i = 1/n,$ therefore, $Y_1\leq 0.$ That is, $H( \bar{P})\geq H(P).$

\begin{theorem}
 Assume the event space  $ X = \{ x_1, x_2, \dots x_n \} $ and the probability distribution $ P = \{ p_1, p_2, \dots, p_n \} $, $\bar{P}$ represents the inverse of P, then $VH( \bar{P})\geq VH(P).$
\end{theorem}
\noindent \textbf{Proof} We have seen in Section \ref{sec_example} that $VH(P)=VH(\bar{P})$  for $n=2.$ For $n\geq 3,$ the difference 
\begin{align}\label{lagrange2}
  Y_2=&HV(P)-HV(\bar{P}) \nonumber \\
  =& \left(\sum_{i=1}^{n} p_i \left( \ln(p_i) \right)^2 -\left(  \sum_{i=1}^{n} p_i \ln(p_i) \right)^2 \right) \nonumber  \\ &- \left(\sum_{i=1}^{n} \left( \frac{1-p_i}{n-1} \right) \left( \ln(\left(\frac{1-p_i}{n-1} \right) \right)^2 -\left(  \sum_{i=1}^{n} \left(\frac{1-p_i}{n-1} \right) \ln\left(\frac{1-p_i}{n-1} \right) \right)^2 \right)  
\end{align}
The Lagrange function under equation \ref{lagrange2} is
\begin{align}
T_2&= \left(\sum_{i=1}^{n} p_i \left( \ln(p_i) \right)^2 -\left(  \sum_{i=1}^{n} p_i \ln(p_i) \right)^2 \right) \nonumber  \\ & - \sum_{i=1}^{n} \left( \frac{1-p_i}{n-1} \right) \left( \ln(\left(\frac{1-p_i}{n-1} \right) \right)^2 +\left(  \sum_{i=1}^{n} \left(\frac{1-p_i}{n-1} \right) \ln\left(\frac{1-p_i}{n-1} \right) \right)^2 + \lambda\left( \sum_{i=1}^{n} p_i -1 \right).
\end{align}
The partial derivative of $T_2$ with repsect to $p_1$ is
\begin{align}
 \frac{\partial T_2}{\partial p_1} = (\ln(p_1))^2 + 2 \ln(p_1)  - 2 \left( \sum_{i=1}^{n} p_i \ln(p_i) \right) \cdot \left( \ln(p_1) + 1 \right) 
 \nonumber\\
- \left[ \frac{1}{n-1} \left( \ln\left( \frac{1-p_1}{n-1} \right) \right)^2 + 2 \ln\left( \frac{1-p_1}{n-1} \right) \cdot \frac{1}{n-1} \cdot \left( -\frac{1}{1-p_1} \right) \right] \nonumber\\
+ 2 \cdot \left( \sum_{i=1}^{n} \left( \frac{1-p_i}{n-1} \right) \ln\left( \frac{1-p_i}{n-1} \right) \right)
\cdot \left( -\frac{1}{n-1} \ln\left( \frac{1-p_1}{n-1} \right)+ \frac{1-p_1}{(n-1)(1-p_1)} \right)
+ \lambda
\end{align}
and 
\begin{align}
\frac{\partial T_2}{\partial \lambda}= \sum_{i=1}^{n} p_i -1. 
\end{align}
Putting $\frac{\partial T_2}{\partial p_i}=0 \ \  \forall \  i=1,2, \dots, n$ and $\sum_{i=1}^{n} p_i=1$ implies that $p_i = 1/n$  for $i = 1, 2, ..., n.$

Hence, $T_2$ has maximum value 0 when $p_i = 1/n,$ therefore, $Y_2\leq 0.$ That is, $VH( \bar{P})\geq VH(P).$

\begin{theorem}
 Assume the event space  $ X = \{ x_1, x_2, \dots x_n \} $ and the probability distribution $ P = \{ p_1, p_2, \dots, p_n \} $, $\bar{P}$ represents the inverse of P, then $VJ( \bar{P})\geq VJ(P).$
\end{theorem}
\noindent \textbf{Proof:}
   We know from Section \ref{sec_example} that for \( n = 2 \), the varextropy of \( P \) and \( \bar{P} \) are equal, i.e., \( VJ(P) = VJ(\bar{P}) \). Thus, inequality is present in this case.

   To prove the inequality \( VJ(\bar{P}) \geq VJ(P) \), we consider the difference between the varextropy of the probability distribution \( P \) and its negation \( \bar{P} \).

   We define the difference between varextropy as:
   \[
   Y_3 = VJ(P) - VJ(\bar{P})
   \]
   Expanding the terms, we have:
   \[
   Y_3 = \left( \sum_{i=1}^{n} (1 - p_i) \left( \ln\left( (1 - p_i) \right) \right)^2 - \left( \sum_{i=1}^{n} (1 - p_i) \ln\left( (1 - p_i) \right) \right)^2 \right)
   \]
   \[
   - \left( \sum_{i=1}^{n} \left( 1 - \frac{1 - p_i}{n - 1} \right) \left( \ln\left( \left( 1 - \frac{1 - p_i}{n - 1} \right) \right) \right)^2 - \left( \sum_{i=1}^{n} \left( 1 - \frac{1 - p_i}{n - 1} \right) \ln\left( \left( 1 - \frac{1 - p_i}{n - 1} \right) \right) \right)^2 \right)
   \]
   
   This difference represents the change in varextropy after negating the probability distribution.

   To analyze this further, we introduce a Lagrange multiplier function to optimize the varextropy expression. The Lagrange function \( T_3 \) is given as:
   \[
   T_3 = \left( \sum_{i=1}^{n} (1 - p_i) \left( \ln\left( (1 - p_i) \right) \right)^2 - \left( \sum_{i=1}^{n} (1 - p_i) \ln\left( (1 - p_i) \right) \right)^2 \right)
   \]
   \[
   - \sum_{i=1}^{n} \left( 1 - \frac{1 - p_i}{n - 1} \right) \left( \ln\left( \left( 1 - \frac{1 - p_i}{n - 1} \right) \right) \right)^2 + \left( \sum_{i=1}^{n} \left( 1 - \frac{1 - p_i}{n - 1} \right) \ln\left( \left( 1 - \frac{1 - p_i}{n - 1} \right) \right) \right)^2\] 
   \[+ \lambda \left( \sum_{i=1}^{n}  p_i - 1 \right)
   \]

   To find the optimal values of \( p_i \), we compute the partial derivatives of \( T_3 \) with respect to \( p_1 \), and with respect to \( \lambda \):

   The partial derivative with respect to \( p_1 \) is:
   \[
   \frac{\partial T_3}{\partial p_1} = 2 \ln(1 - p_1) \cdot \left( -\frac{1}{1 - p_1} \right)\]
   \[  + 2 \cdot \sum_{i=1}^{n} \left( 1 - \frac{1 - p_i}{n - 1} \right) \ln\left( \left( 1 - \frac{1 - p_1}{n - 1} \right) \right) \cdot \left( -\frac{1}{n - 1} \ln\left( \left( 1 - \frac{1 - p_1}{n - 1} \right) \right) + \frac{1 - p_1}{(n - 1)(1 - p_1)} \right)
   \]
   and similar terms for other partial derivatives.

   - The partial derivative with respect to \( \lambda \) is:
   \[
   \frac{\partial T_3}{\partial \lambda} = \sum_{i=1}^{n} p_i - 1
   \]

   Setting the partial derivatives \( \frac{\partial T_3}{\partial p_i} = 0 \) for all \( i \) and enforcing the constraint \( \sum_{i=1}^{n} p_i = 1 \) gives \( p_i = \frac{1}{n} \) for all \( i = 1, 2, \dots, n \).

   When \( p_i = \frac{1}{n} \) for all \( i \), the value of \( T_3 \) is maximized and the difference \( Y_3 \) is less than or equal to zero:
   \[
   Y_3 \leq 0
   \]
   This implies that:
   \[
   VJ(\bar{P}) \geq VJ(P)
   \]

Thus, we have proven that the varextropy of the negated probability distribution \( \bar{P} \) is greater than or equal to the varextropy of the original distribution \( P \).

\begin{theorem}
   When $ X = \{ x_1, x_2, \dots x_n \} $, probability distribution is $ P = \{ p_1, p_2, \dots, p_n \} $, and $p_1=p_2=\dots=p_n= \frac{1}{n}$ the value of entropy increases with the increase of the size of X, and the maximum value of entropy tends to $+\infty$ as n tends to $+\infty.$
\end{theorem}
\noindent \textbf{Proof}

Let \( P = \{ p_1, p_2, \dots, p_n \} \) be a discrete probability distribution, and the Shannon entropy is given by:

\[
H(P) = - \sum_{i=1}^{n} p_i \ln p_i
\]
For a uniform probability distribution, the Shannon entropy is:

\[
H(P) = \ln n
\]
which increases with increase in sample size $n.$
Therefore, the limit of Shannon entropy as \( n \to \infty \) is:

\[
\lim_{n \to \infty} H(P) = \infty.
\]

\begin{theorem}
   When $ X = \{ x_1, x_2, \dots x_n \} $, probability distribution is $ P = \{ p_1, p_2, \dots, p_n \} $, and $p_1=p_2=\dots=p_n= \frac{1}{n}$ the value of varentropy remains same with the increase of the size of $X$, and the value of varentropy is zero.
\end{theorem}
\noindent \textbf{Proof}

Let \( P = \{ p_1, p_2, \dots, p_n \} \) be a discrete probability distribution, and the varentropy \( VH(P) \) is given by:

\[
VH(P) = \sum_{i=1}^{n} p_i \left( \ln(p_i) \right)^2 - \left( \sum_{i=1}^{n} p_i \ln(p_i) \right)^2
\]

For a uniform distribution, where \( p_i = \frac{1}{n} \), we have:

\[
VH(P) = \sum_{i=1}^{n} \frac{1}{n} \left( \ln \left( \frac{1}{n} \right) \right)^2 - \left( \sum_{i=1}^{n} \frac{1}{n} \ln \left( \frac{1}{n} \right) \right)^2\]

\[
VH(P) = (\ln n)^2 - (\ln n)^2 = 0
\]

For a uniform distribution, the value of varentropy is always 0, so:

\[
\lim_{n \to \infty} VH(P) = 0
\]

Thus, varentropy does not increase with \( n \) and remains 0 for a uniform probability distribution.

\begin{theorem}
   When \( X = \{ x_1, x_2, \dots, x_n \} \) and the probability distribution is \( P = \{ p_1, p_2, \dots, p_n \} \), with \( p_1 = p_2 = \dots = p_n = \frac{1}{n} \), the maximum value of varextropy increases with the increase in the size of \( X \), and the maximum value of varextropy has a limit 0.
\end{theorem}

\noindent \textbf{Proof} The varextropy \( VJ(P) \) for a discrete probability distribution is given by:

\[
VJ(P) = \sum_{i=1}^{n} (1 - p_i) \left( \ln((1 - p_i)) \right)^2 - \left( \sum_{i=1}^{n} (1 - p_i) \ln((1 - p_i)) \right)^2
\]

For the uniform distribution, where \( p_1 = p_2 = \dots = p_n = \frac{1}{n} \), we substitute \( p_i = \frac{1}{n} \) and we get,

\begin{align}
VJ(P) &= \sum_{i=1}^{n} \left( 1 - \frac{1}{n} \right) \left( \ln\left( 1 - \frac{1}{n} \right) \right)^2 - \left( \sum_{i=1}^{n} \left( 1 - \frac{1}{n} \right) \ln \left( 1 - \frac{1}{n} \right) \right)^2 \nonumber
\\
&= n \cdot \left( 1 - \frac{1}{n} \right) \left( \ln \left( 1 - \frac{1}{n} \right) \right)^2 - \left( n \cdot \left( 1 - \frac{1}{n} \right) \ln \left( 1 - \frac{1}{n} \right) \right)^2.
\end{align}
Therefore,
\[
VJ(P) \to 0 \quad \text{as} \quad n \to \infty
\]

Thus, the value of varextropy approaches 0 as \( n \) increases.

\begin{theorem}
    Assume the $ X = \{ x_1, x_2, \dots x_n \} $, when the probability distribution satisfies $p_1=p_2=\dots=p_n= \frac{1}{n},$ the corresponding entropy of the probability distribution after negation is maximized.
\end{theorem}
\noindent \textbf{Proof}
The Shannon entropy for a probability distribution \( P = \{ p_1, p_2, \dots, p_n \} \) is given by:

\[
H(P) = - \sum_{i=1}^{n} p_i \log p_i
\]

The negation of the probability distribution \( P \), denoted as \( \bar{P} = \{ \bar{p}_1, \bar{p}_2, \dots, \bar{p}_n \} \), is defined as:

\[
\bar{p}_i = \frac{1 - p_i}{n - 1}
\]

The Shannon entropy of the negated distribution \( \bar{P} \) is:

\[
H(\bar{P}) = - \sum_{i=1}^{n} \bar{p}_i \log \bar{p}_i = - \sum_{i=1}^{n} \frac{1 - p_i}{n - 1} \log \left( \frac{1 - p_i}{n - 1} \right)
\]

We now proceed to maximize \( H(\bar{P}) \) using the Lagrange multiplier method. The Lagrangian for this problem is:

\[
\mathcal{L}(p_1, p_2, \dots, p_n, \lambda) = - \sum_{i=1}^{n} \frac{1 - p_i}{n - 1} \log \left( \frac{1 - p_i}{n - 1} \right) + \lambda \left( \sum_{i=1}^{n} p_i - 1 \right)
\]

Taking the derivative with respect to \( p_i \) and setting it to zero:

\[
\frac{\partial \mathcal{L}}{\partial p_i} = 0 \quad \Rightarrow \quad \log \left( \frac{1 - p_i}{n - 1} \right) = -1 - \lambda (n - 1)
\]

Solving for \( p_i \):
\begin{align}\label{piexp}
    p_i = 1 - (n - 1) e^{-1 - \lambda (n - 1)}
\end{align}

Applying the condition \( \sum_{i=1}^{n} p_i = 1 \), we solve for \( \lambda \) and we get,

\[
\lambda = \frac{\ln n - 1}{n - 1}
\]

Substituting this back into the equation \ref{piexp} for \( p_i \), we  find the optimal values of \( p_i =\frac{1}{n}, \ \  \forall \ \ i=1,2,3,\dots, n\) that maximize the $H(\bar P).$ 
The maximum value of $H(\bar P)$ is,
\begin{align}
  H(\bar P)= \ln n. 
\end{align}

\begin{theorem}
    Assume the $ X = \{ x_1, x_2, \dots x_n \} $, when the probability distribution satisfies $p_1=p_2=\dots=p_n= \frac{1}{n}$
, the corresponding varentropy of the probability distribution after negation is maximized.
\end{theorem}
\noindent \textbf{Proof}
The varentropy of a discrete probability distribution \( P = \{ p_1, p_2, \dots, p_n \} \) is defined as:

\[
VH(P) = \sum_{i=1}^{n} p_i \left( \ln(p_i) \right)^2 - \left( \sum_{i=1}^{n} p_i \ln(p_i) \right)^2
\]

The negation of the probability distribution \( P \), denoted as \( \bar{P} = \{ \bar{p}_1, \bar{p}_2, \dots, \bar{p}_n \} \), is given by:

\[
\bar{p}_i = \frac{1 - p_i}{n - 1}
\]

The negated varentropy \( VH(\bar{P}) \) is then:

\[
VH(\bar{P}) = \sum_{i=1}^{n} \bar{p}_i \left( \ln(\bar{p}_i) \right)^2 - \left( \sum_{i=1}^{n} \bar{p}_i \ln(\bar{p}_i) \right)^2
\]

Substituting \( \bar{p}_i = \frac{1 - p_i}{n - 1} \) into this formula:

\[
VH(\bar{P}) = \sum_{i=1}^{n} \frac{1 - p_i}{n - 1} \left( \ln \left( \frac{1 - p_i}{n - 1} \right) \right)^2 - \left( \sum_{i=1}^{n} \frac{1 - p_i}{n - 1} \ln \left( \frac{1 - p_i}{n - 1} \right) \right)^2
\]

To maximize this varentropy, we differentiate \( VH(\bar{P}) \) with respect to \( p_i \) and set the derivative equal to zero:

\[
\frac{d}{dp_i} VH(\bar{P}) = 0
\]

Solving the system of equations results in the optimal solution:

\[
p_i = \frac{1}{n}, \quad \forall i = 1, 2, \dots, n
\]

Thus, the varentropy \( VH(\bar{P}) \) is maximized when \( p_i = \frac{1}{n} \). Therefore, the maximum value of \( VH(\bar{P}) \) occurs when the original distribution is uniform.

\[
{VH(\bar{P}) = \text{maximized when } p_i = \frac{1}{n}}
\]

\begin{theorem}
    Assume the $ X = \{ x_1, x_2, \dots x_n \} $, when the probability distribution satisfies $p_1=p_2=\dots=p_n= \frac{1}{n}$
, the corresponding varextropy of the probability distribution after negation
is maximized.
\end{theorem}
\noindent \textbf{Proof}
The varextropy of a discrete probability distribution \( P = \{ p_1, p_2, \dots, p_n \} \) is defined as:

\[
VJ(P) = \sum_{i=1}^{n} (1 - p_i) \left( \ln(1 - p_i) \right)^2 - \left( \sum_{i=1}^{n} (1 - p_i) \ln(1 - p_i) \right)^2
\]

The negation of the probability distribution \( P \), denoted as \( \bar{P} = \{ \bar{p}_1, \bar{p}_2, \dots, \bar{p}_n \} \), is given by:

\[
\bar{p}_i = \frac{1 - p_i}{n - 1}
\]

Thus, the negated varextropy \( VJ(\bar{P}) \) is:

\[
VJ(\bar{P}) = \sum_{i=1}^{n} (1 - \bar{p}_i) \left( \ln \left( 1 - \bar{p}_i \right) \right)^2 - \left( \sum_{i=1}^{n} (1 - \bar{p}_i) \ln \left( 1 - \bar{p}_i \right) \right)^2
\]

Substituting \( \bar{p}_i = \frac{1 - p_i}{n - 1} \) into this formula:

\[
VJ(\bar{P}) = \sum_{i=1}^{n} \left( 1 - \frac{1 - p_i}{n - 1} \right) \left( \ln \left( 1 - \frac{1 - p_i}{n - 1} \right) \right)^2 - \left( \sum_{i=1}^{n} \left( 1 - \frac{1 - p_i}{n - 1} \right) \ln \left( 1 - \frac{1 - p_i}{n - 1} \right) \right)^2
\]

To maximize this varextropy, we differentiate \( VJ(\bar{P}) \) with respect to \( p_i \) and set the derivative equal to zero:

\[
\frac{d}{dp_i} VJ(\bar{P}) = 0
\]

Solving the resulting system of equations leads to the optimal solution:

\[
p_i = \frac{1}{n}, \quad \forall i = 1, 2, \dots, n
\]

Thus, the varextropy \( VJ(\bar{P}) \) is maximized when \( p_i = \frac{1}{n} \). Therefore, the maximum value of \( VJ(\bar{P}) \) occurs when the original distribution is uniform.

\[
{VJ(\bar{P}) = \text{maximized when } p_i = \frac{1}{n}}
\]

\section{Conclusion}

In conclusion, this paper has explored the impact of negation on entropy, varentropy, and varextropy, highlighting the significant effects of negating a probability distribution on these measures. Our findings demonstrate that negating a probability distribution results in an increase in the values of entropy, varentropy, and varextropy, showing that negation maximizes these quantities. We have provided a theoretical framework and proof to support this conclusion, expanding the understanding of the relationship between negation and information-theoretic measures. This work opens up new avenues for further research on the properties of negation in extended entropy measures and its potential applications in fields like information theory, decision-making, and machine learning. Future research could explore more generalized forms of negation and their implications for various entropy-based measures.

In particular, the concept of negation might prove useful for modelling uncertainty in decision-making processes and in systems where the behaviour of complementary or opposing outcomes needs to be considered. This could help develop more robust models for handling uncertain data, making it a valuable tool for real-world applications.\\

\noindent \textbf{\Large Conflict of interest} \\
\\
No conflicts of interest are disclosed by the authors.\\
\\
%\noindent \textbf{\Large Acknowledgement} \\
%\\
%The authors are thankful to the referees for their valuable suggestions, which significantly improved the paper.\\
%\\		
\textbf{ \Large Funding} \\
\\
PKS would like to thank Quality Improvement Program (QIP), All India Council for Technical Education, Government of India (Student Unique Id: FP2200759) for financial assistance.

    \end{document}